\documentclass[12pt,nosumlimits,nonamelimits]{amsart}
\usepackage{mathtools}

\newtheorem*{thm*}{Theorem}
\newcommand{\IB}{\mathbb{B}}

\newcommand{\IF}{\mathbb{F}}
\newcommand{\IM}{\mathbb{M}}

\newcommand{\IP}{\mathbb{P}}
\newcommand{\IZ}{\mathbb{Z}}

\newcommand{\cH}{\mathcal{H}}
\newcommand{\cK}{\mathcal{K}}
\newcommand{\cL}{\mathcal{L}}
\newcommand{\cM}{\mathcal{M}}
\newcommand{\cR}{\mathcal{R}}

\newcommand{\id}{\mathrm{id}}
\DeclareMathOperator{\Tr}{Tr}
\DeclareMathOperator{\SL}{SL}

\title[Non-quasidiagonality of $\mathcal{R}$]{The hyperfinite $\mathrm{II}_1$ factor is not quasidiagonal}
\author{Narutaka Ozawa}
\address{RIMS, Kyoto University, \mbox{606-8502} Japan}
\email{narutaka@kurims.kyoto-u.ac.jp}
\thanks{The author was partially supported by JSPS KAKENHI Grant Numbers 24K00527, 25H00588, 25H00593}

\date{\today}

\begin{document}
\begin{abstract}
We construct an MF $\mathrm{C}^*$-alge\-bra $A$ such that the spatial tensor 
product $A\otimes \cR$ of $A$ and the hyperfinite $\mathrm{II}_1$ factor $\cR$ 
contains a proper isometry. Consequently, stable finiteness of 
$\mathrm{C}^*$-alge\-bras is not stable under tensor product and 
the hyperfinite $\mathrm{II}_1$ factor is not a quasidiagonal $\mathrm{C}^*$-alge\-bra. 
The $\mathrm{C}^*$-alge\-bra $A$ is non-simple and 
has only non-faithful tracial states. 
This result was obtained using OpenAI’s Chat GPT Pro 6.0.
\end{abstract}
\maketitle
Let $A$ be a $\mathrm{C}^*$-alge\-bra, 
which is assume to be unital and separable for a simpler presentation. 
Let $\cM \coloneq \prod_n \IM_n$ be the $\ell_\infty$-direct sum of matrix algebras 
and $\cK \coloneq \bigoplus_n \IM_n$ be the ideal of the $c_0$-sum. 
We denote by $Q\colon\cM\to\cM/\cK$ the quotient map. 
The $\mathrm{C}^*$-alge\-bra $A$ is said to 
be \emph{quasidiagonal} (Exercise 7.1.3 in \cite{bo}) 
if there is a contractive completely positive map 
$\Theta\colon A \to\cM$ such that $Q\circ\Theta$ 
is a (not necessarily unital) faithful $*$-homo\-mor\-phism. 
(A not necessarily separable $\mathrm{C}^*$-alge\-bra, 
such as the hyperfinite $\mathrm{II}_1$ factor $\cR$, is 
quasidiagonal if every separable $\mathrm{C}^*$-sub\-alge\-bra 
of it is quasidiagonal.)
Every quasidiagonal $\mathrm{C}^*$-alge\-bra $A$ is 
\emph{MF} (Definition 3.2.1 in \cite{bk}) in the sense that 
it is $*$-iso\-mor\-phic to a (possibly non-unital) 
$\mathrm{C}^*$-sub\-alge\-bra of $\cM/\cK$. 
The converse need not be true (\cite{wassermann}). 
A $\mathrm{C}^*$-alge\-bra is said to be \emph{infinite} 
if it contains a proper isometry; else it is said to be \emph{finite}. 
Every MF algebra $A$ is \emph{stably finite}, i.e., 
$\IM_n \otimes A$ is finite for every $n$. 
The converse need not be true (\cite{jn+}, 
see the remark after Problem VII in \cite{stw}). 
The problem whether stable finiteness 
of \emph{simple} $\mathrm{C}^*$-alge\-bras 
is stable under the spatial tensor product has attracted considerable 
attention (see e.g., \cite{bk,gs,mr,stw}) in connection with Kaplansky's problem 
on $\mathrm{AW}^*$ factors. 
Here we answer in the negative the less interesting analogue of this problem 
for a \emph{non-simple} $\mathrm{C}^*$-alge\-bra and 
the problem in Section 6.6 in \cite{brown} (see also 
Problem 10.4.9 in \cite{bo} and Problem X in \cite{stw}).
In fact our example has only non-faithful tracial states, 
and thus does not resolve Kaplansky's problem. 

\begin{thm*}
There is a unital separable MF $\mathrm{C}^*$-alge\-bra $A$ 
such that the spatial tensor product $A\otimes \cR$ 
of $A$ and the hyperfinite $\mathrm{II}_1$ factor $\cR$ 
is infinite. Consequently, stable finiteness of 
$\mathrm{C}^*$-alge\-bras is not stable under tensor product and 
the hyperfinite $\mathrm{II}_1$ factor is not a quasidiagonal $\mathrm{C}^*$-alge\-bra. 
\end{thm*}

We start the proof of Theorem. 
Our construction builds upon the developments of \cite{wassermann}.  
For a prime number $p$, the standard action of 
the group $G \coloneq \SL(3,\IZ)$ on the finite projective plane 
$\IP^2(\IF_p) \coloneq (\IF_p^3\setminus\{0\})/\IF_p^\times$ 
is doubly transitive. Hence the corresponding unitary 
representation $\pi$ of $G$ 
on $\cH \coloneq \ell_2(\IP^2(\IF_p))\cap\{\mathbf{1}\}^\perp$ is irreducible 
with $\dim\cH = p^2+p$. 
We denote by $p_k$ the $k$-th prime number and by $(\pi_k,\cH_k)$ 
the above unitary representation of $G$ associated with $p_k$. 
We denote the trivial one-dimensional representation by $(\pi_0,\cH_0)$. 
All we need below is that $G$ has Kazhdan's property (T), $\pi_k$ are irreducible, and 
$d_k \coloneq \dim \cH_k$ satisfies $d_k\nearrow\infty$ 
and $\sup d_k/d_{k-1} < \infty$ (by Bertrand's postulate). 

We consider the complex conjugate representations $\bar{\pi}_k$ on $\bar{\cH}_k$, 
which is same as $(\pi_k,\cH_k)$ in our setting of orthogonal representations. 
By Schur's lemma $(\pi_k\otimes\bar{\pi}_l)(G)$, $k\neq l$ does not have 
nonzero invariant vectors and $(\pi_k\otimes\bar{\pi}_k)(G)$ has unique 
invariant vector ``the identity vector'' (as a Hilbert--Schmidt operator) 
$\sum_i \zeta_i\otimes\bar{\zeta}_i$, up to the scalar multiple, 
where $(\zeta_i)_i$ is any orthonormal basis. 
The rank one projection $p_{\cH} \in \IB(\cH\otimes\bar{\cH})$ 
corresponding to the identity vector 
satisfies $p_{\cH} = (\dim\cH)^{-1} \sum_{i,j} e_{i,j} \otimes \bar{e}_{i,j}$ 
and $(\id\otimes\Tr)(p_{\cH})=(\dim\cH)^{-1}$ in $\IB(\cH)$. 
Here $\{e_{i,j}\}$ is any matrix unit for $\IB(\cH)$ 
and $\Tr$ denotes the unnormalized trace. 
For every isometry $T\colon \cK \to \cH$ the operator 
$(\dim\cH/\dim\cK)^{1/2} p_{\cH}(T\otimes\bar{T})$ 
is a rank-one partial isometry from 
$\cK\otimes\bar{\cK}$ into $\cH\otimes\bar{\cH}$ 
that intertwines $p_{\cK}$ and $p_{\cH}$. 

Fix a positive integer $c > \sup d_k/d_{k-1}$ 
and a unital normal embedding $\iota$ of $\prod_k\IB(\bar{\cH}_k)$ into $\cR$, 
for which the minimal projections 
in $M_k\coloneq\iota(\IB(\bar{\cH}_k))$ have trace $c^{-k}$, 
where the trace $\tilde{\tau}$ on $\cR$ is normalized 
so that $\tilde{\tau}(1) = \sum_{k=0}^\infty d_k c^{-k} < \infty$.
We view $M_k$ as a non-unital subalgebra of $\cR$ and 
put $\iota_k\coloneq\iota|_{\IB(\bar{\cH}_k)}\colon \IB(\bar{\cH}_k)\to M_k$. 

We denote by $\ell_2(m)$ the $m$-dimensional Hilbert space. 
We consider the Hilbert space 
$\cL_n \coloneq \bigoplus_{k=0}^n \ell_2(c^k) \otimes \cH_k$
and the $\mathrm{C}^*$-alge\-bras 
$\cM \coloneq \prod_n \IB(\cL_n)$ and $\cK \coloneq \bigoplus_n \IB(\cL_n)$. 
The renewal of $\cM$ and $\cK$ does not affect the argument. 
For a sequence $x_n \in \IB(\cL_n)$ with $\sup\|x_n\|<\infty$, 
we will write the corresponding element in $\cM$ as $(x_n)_n$ or $\sum_n x_n$ 
according to our convenience. 
We set 
$\sigma_n(g) \coloneq \bigoplus_{k=0}^n (1\otimes\pi_k)(g) \in \IB(\cL_n)$, 
$\sigma(g) \coloneq (\sigma_n(g))_n \in \cM$, 
and $\bar{\pi}(g) \coloneq \iota( (\bar{\pi}_n(g))_n ) \in \cR$, 
for $g\in G$. 
We fix a finite symmetric generating subset $E\subset G$ that contains $1$. 
Since $G$ has property (T), the self-adjoint element 
\[
\frac{1}{|E|}\sum_{g\in E} \sigma(g) \otimes \bar{\pi}(g) \in \cM\otimes \cR
 \mbox{ acting on }\bigoplus_n \bigoplus_{k=0}^n \ell_2(c^k)\otimes \cH_k\otimes L^2(\cR)
\]
has $1$ as an isolated point of its spectrum (see e.g., Lemma 12.1.8 in \cite{bo}) 
and the corresponding projection $q$ belongs to $\cM\otimes \cR$. 

The projection $q$ decomposes into orthogonal sums 
$q = \sum_n q_n$ and $q_n = \sum_{k=0}^n q_{n,k} \in \IB(\cL_n) \otimes \cR$, 
where $q_{n,k} \in \IB(\ell_2(c^k) \otimes \cH_k) \otimes M_k$ 
is the $c^k$ diagonal sum of the rank one projection 
onto subspace spanned by the identity element in $\cH_k\otimes\bar{\cH}_k$. 
That is, for the orthogonal projection $r_{n,k}\in\IB(\cL_n)$ 
onto $\ell_2(c^k) \otimes \cH_k \subset\cL_n$, one has $q_{n,k} = (r_{n,k} \otimes 1)q$ 
in $\cM\otimes\cR$. 
Note that $q_{n,k}$, $k=0,\ldots,n$, are equivalent 
in $\IB(\cL_n) \otimes \cR$, because 
\[
(\id\otimes\tilde{\tau})(q_{n,k}) = \frac{1}{d_k c^k} r_{n,k}
\ \mbox{ and }\ 
(\Tr\otimes\tilde{\tau})(q_{n,k}) = 1. 
\]
It is clear that the projections 
\[
q_{\mathrm{b}} \coloneq (q_{n,0})_n = ( (r_{n,0})_n \otimes 1) q
\ \mbox{ and }\ 
q_{\mathrm{t}} \coloneq (q_{n,n})_n = ( (r_{n,n})_n \otimes 1) q
\]
belong to $\cM\otimes\cR$. 

We claim that $(Q\otimes\id_{\cR})( q_{\mathrm{b}} ) \neq 0$ 
and $(Q\otimes\id_{\cR})( q_{\mathrm{t}} ) = 0$ in $(\cM/\cK)\otimes\cR$. 
Recall that an element $x\in(\cM\otimes\cR)_+$ 
belongs to $\ker(Q\otimes\id_{\cR})$ 
if and only if $(\id_{\cM}\otimes\tilde{\tau})(x) \in \ker Q$ 
(see e.g., Lemma 4.1.8 in \cite{bo}). 
Thus the claim follows from the following. 
\[
(\id_{\cM}\otimes\tilde{\tau})(q_{\mathrm{b}}) = (r_{n,0})_n \notin \cK
\ \mbox{ and }\ 
(\id_{\cM}\otimes\tilde{\tau})(q_{\mathrm{t}}) = (\frac{1}{d_n c^n}r_{n,n})_n \in \cK.
\]

For the proof of Theorem, it is left to show that 
$1-q_{\mathrm{b}}$ and $1-q_{\mathrm{t}}$ 
are equivalent in $\cM\otimes\cR$ (not only in the von Neumann 
algebra $\cM\mathbin{\bar{\otimes}}\cR$). 
Then, non-quasidiagonality of $\cR$ follows from Lemma~4 in \cite{gs}.
We view $\ell_2(c^k)=\ell_2(c)\otimes\ell_2(c^{k-1})$ 
and consider the isometries 
$s_{k,i}\colon \ell_2(c^{k-1}) \to \ell_2(c)\otimes\ell_2(c^{k-1})$, 
$\xi\mapsto \delta_i\otimes\xi$. 
We also fix isometries $t_k\colon\cH_{k-1}\to\cH_k$. 
We define $a^{(n)}_i\in\IB(\cL_n)$ by 
\[
a^{(n)}_i \coloneq \sum_{k=1}^n (\frac{d_k}{d_{k-1}})^{1/2} s^{(n)}_{k,i} \otimes t^{(n)}_k.
\]
Here $s^{(n)}_{k,i}$ and $t^{(n)}_k$ are the copies 
of $s_{k,i}$ and $t_k$, viewed as operators acting on 
the Hilbert spaces appearing in the decomposition 
$\cL_n = \bigoplus_{k=0}^n \ell_2(c^k) \otimes \cH_k$.
Thus the operator $a^{(n)}_i$ shifts the direct summand of 
$\cL_n$ forward by $k-1\to k$.
Put $a_i \coloneq (a^{(n)}_i)_n \in\cM$ by noticing that $\sup_n\| a^{(n)}_i \|^2 < c$. 
Let $k\geq 1$ and $z_{k-1}$ 
denote the unit for $M_{k-1}$ in $\prod_k M_k\subset\cR$. 
We fix equivalent and mutually orthogonal family $(f_{k-1,i})_{i=1}^c$ of projections 
in $M_{k-1}'\cap z_{k-1}\cR z_{k-1}$ with sum $z_{k-1}$. The minimal projections in $f_{k-1,i}M_{k-1}$ 
and the minimal projections in $M_k$ have the same trace values 
(which is $c^{-k}$ in our normalization). 
Thus there is a partial isometry $b_{k,i}$ in $\cR$ such that 
$b_{k,i}^*b_{k,i} = z_{k-1}f_{k-1,i}$, 
$b_{k,i}b_{k,i}^* = \iota_k(\bar{t}_k\bar{t}_k^*) =: z_k' \le z_k$, and 
$b_{k,i}^*\iota_k(\bar{x})b_{k,i} = \iota_{k-1}(\bar{t}_k^* \bar{x} \bar{t}_k)f_{k-1,i}$ 
for every $\bar{x}\in\IB(\bar{\cH}_k)$. 
Put $b_i \coloneq \sum_{k=1}^\infty b_{k,i} \in \cR$ 
and $w\coloneq \sum_{i=1}^c a_i \otimes b_i \in \cM\otimes\cR$. 
Since $b_{k,i}$, $k=1,2,\ldots$, are partial isometries in $\cR$ 
with mutually orthogonal initial spaces and mutually orthogonal 
ranges, the infinite sum $b_i$ makes sense in $\cR$ 
as the limit in the strong operator topology. 
We claim that $qw$ is a partial isometry 
such that $(qw)^*(qw) = q - q_{\mathrm{t}}$ 
and $(qw)(qw)^* = q - q_{\mathrm{b}}$.
Indeed, 
\[
q_{n,k} w
 = q_{n,k} \sum_i a^{(n)}_i \otimes b_{k,i} 
 = \sum_i q_{n,k} w_{n,k,i}
\]
where $w_{n,k,i} = (d_k/d_{k-1})^{1/2} s^{(n)}_{k,i} \otimes t^{(n)}_k \otimes b_{k,i}$ for $k\geq1$. 
Note that $q_{n,0}w=0$. 
Moreover, by the previous discussion on the rank-one projections 
associated with the identity vectors, one has 
\[
(q_{n,k} w_{n,k,i})^*(q_{n,k} w_{n,k,i}) = q_{n,k-1}(1\otimes f_{k-1,i})
\]
and 
\[
(q_{n,k} w_{n,k,i})(q_{n,k} w_{n,k,i})^* = (s^{(n)}_{k,i})(s^{(n)}_{k,i})^* q_{n,k}.
\]
Hence 
\[
(q_n w)^*(q_n w) =\sum_{k=1}^n (q_{n,k}w)^*(q_{n,k}w)
 = \sum_{k=1}^n q_{n,k-1} = q_n - q_{n,n}
\]
and 
\[
(q_n w)(q_n w)^* = \sum_{k=1}^n (q_{n,k}w)(q_{n,k}w)^*
 = \sum_{k=1}^n q_{n,k} = q_n-q_{n,0}.
\]
This proves the claim. 
It follows that $1-q_{\mathrm{t}}$ and $1-q_{\mathrm{b}}$ 
are equivalent via the partial isometry $(1-q)+qw$. 
Finally, we define the separable MF algebra $A$ 
by $A \coloneq B/\cK$, where $B$ is the 
$\mathrm{C}^*$-sub\-alge\-bra of $\cM$ 
generated by $\sigma(G)$, $\{ a_i \}_{i=1}^c$, and $\cK$. 
\hspace*{\fill}\qedsymbol
\smallskip

We note that a similar but simpler proof shows that 
the spatial tensor product $A\otimes \bar{A}$ of 
an MF algebra $A$ with its opposite is also infinite. 
See Remark 3.7 in \cite{rs}. 
Here is a sketch of the proof. 
Since we do not work with $\cR$, 
forget about $\ell_2(c^k)$ and consider 
$\cL_n \coloneq \bigoplus_{k=0}^n \cH_k$. 
Thus $\cM\otimes\bar{\cM}$ naturally acts 
on $\bigoplus_{n,m}\cL_n\otimes\bar{\cL}_m$. 
Put $a \coloneq (\sum_{k=1}^n (\frac{d_k}{d_{k-1}})^{1/4} t^{(n)}_k)_n \in \cM$, 
where $t^{(n)}_k$ is the copy of isometry 
$t_k$ from $\cH_{k-1}$ into $\cH_k$, as before. 
Then the spectral projection $p$ corresponding to 
the $(\sigma\otimes\bar{\sigma})(G)$-invariant vectors 
belongs to $\cM\otimes\bar{\cM}$ and 
$w \coloneq p(a\otimes\bar{a}) \in \cM \otimes \bar{\cM}$ 
is a partial isometry 
such that $w^*w =p-p_{\mathrm{t}}$ and $ww^*=p-p_{\mathrm{b}}$, where $p_{\mathrm{t}}$ (resp.\ $p_{\mathrm{b}}$) 
is the projection corresponding to the identity vectors 
in $\cH_n \otimes \bar{\cH}_n$ (resp.\ 
$\cH_0 \otimes \bar{\cH}_0$) that 
sits inside $\cL_n\otimes\bar{\cL}_n$ with $n=0,1,2,\ldots$.
A limit vector state on $\cM$ is a limit-point 
of the vector states associated 
with unit vectors $\xi_n\in\cL_n$ as $n\to\infty$. 
They vanish on $\cK$ and form a faithful family on $\cM/\cK$.
Thus for $x\in(\cM\otimes\bar{\cM})_+$ one has 
$(Q\otimes\bar{Q})(x)=0$ if and only if $(\phi\otimes\bar{\psi})(x)=0$
for all limit vector states $\phi$ and $\psi$. From this, 
it is not hard to see that $(Q\otimes\bar{Q})(p_{\mathrm{t}})=0$ 
and $(Q\otimes\bar{Q})(p_{\mathrm{b}})\neq 0$.

\subsection*{Acknowledgment} 
The outline of the proof of Theorem was conceived 
by the author a long time ago when he tried to solve 
Problem 10.4.9 in \cite{bo}, but he could not bring it to fruition. 
This research was carried out with substantial assistance 
from Chat GPT Pro 6.0 through interactive communication. 
Most of proofs were provided in the end by Chat GPT Pro 6.0.

\end{document}